\documentclass[twoside,12pt]{article}
\usepackage{amsmath}
\usepackage{enumitem}
\usepackage{amssymb}  
\newtheorem{theorem}{Theorem}[section]
\newtheorem{lemma}{Lemma}[section]
\newtheorem{proposition}{Proposition}[section]
\newtheorem{definition}{Definition}[section]

\def\Sum{\displaystyle\sum}

\def\Sup{\displaystyle\sup}
\def\Max{\displaystyle\max}
\def\Inf{\displaystyle\inf}
\def\Lim{\displaystyle\lim}

\DeclareMathOperator*{\esssup}{ess\,sup}
\begin{document}

\begin{center}
{\bf  Analysis of  a free boundary problem on stratified Lie group.   }\vskip 0.5cm {\bf By} \vskip
1cm {\bf Sabri BENSID }\\
\
\\ {\em  Address} : Laboratoire d'Analyse Nonlinéaire et Mathématiques Appliquées,\\ Department of Mathematics,\\ Faculty of Sciences,\\ University of Tlemcen,\\ B.P.119, Tlemcen 13000,
Algeria.\\  Mail: $ edp\_sabri@yahoo.fr$

\end{center}
\begin{abstract}
We present a variational framework for studying the existence of solutions of a class of elliptic free boundary problems on stratified Lie groups. Using the important monotonicity result  in a Non-Euclidean setup, we prove that our solution is   the limit of mountain pass points of a sequence of
$C^1$-functionals approximating the energy .
\end{abstract}

\section{Introduction}
It is well-known that the free boundary problems arise in many models in physics, biology, finance and industry. Such problems can be described by PDE with a priori unknown (free ) interface or boundary. The most important question is to understand the structure and regularity of the free boundary. We let the reader to references \cite{alt-caffarelli,alt-caffarelli-fried,caffa-annal, caff-salsa,Frie-Livre,KINDER} for a brief panorama and for some applications in nature and technology.\\\\
Recently, many authors have attempted to generalize some free boundary problems in the framework of the Heisenberg group like \cite{valdinoc-Ferrari}, when we try to extend the existence results to more general setting of the real sub-Laplacian on stratified Lie group. This paper is a further step in the development of the theory of free boundary value problems in the non-Euclidean setup.\\
More precisely,  we consider a free boundary problem of the type
 \begin{equation}\label{1}
  \left\{\begin{array}{ll}  -\mathcal{L}u=\lambda \chi_{\{u>1\}}g(x,(u-1)_+) & \quad \mbox{
in }\ \Omega\setminus F(u),\\[0.3cm]| \nabla_\mathbb{G}u^+|^2-| \nabla_\mathbb{G}u^-|^2=2 &\quad \mbox{ in }\
  F(u),
\\[0.3cm]u =0 &\quad \mbox{ on }\
  \partial \Omega,
\end{array}
\right.
 \end{equation}
where $\lambda>0,$ $(u-1)_+=\Max(u-1,0),$ $F(u)=\partial \{u>1\}$ is the free boundary of $u$ and $\chi_{\{u>1\}}$ is the characteristic
function of the set.  The domain $\Omega\subset \mathbb{G}$  is a bounded where  $\mathbb{G}$ is a stratified Lie group and $\mathcal{L}$ is the sub-Laplacian which will be defined in the following Section. Finally, $\nabla_\mathbb{G}u^{\pm}$ are the limit of $\nabla_\mathbb{G}u$ from the sets $\{u>1\}$ and $\{u\leq 1\}^{\circ}$ respectively,
and $g:\Omega\times [0,\infty)\rightarrow [0,\infty)$ is a locally H\"{o}lder continuous function satisfying\\\\
$(g_1)$ Fore some $\alpha,\beta>0$ and $1<m<2,$
$$g(x,s)\leq \alpha+\beta s^{m-1},\quad \quad  \forall (x,s)\in \Omega\times[0,\infty).$$
$(g_2)\quad g(x,s)>0,\quad \forall x\in \Omega$ and $s>0.$\\\\
Motivated by the work of Choudhuri and Repovs in \cite{Repsov}, when the authors studied the following problem
\begin{equation}\label{2}
  \left\{\begin{array}{ll}  -\mathcal{L}u=\lambda (u-1)_+^2 f & \quad \mbox{
in }\ \Omega\setminus F(u),\\[0.3cm]| \nabla_\mathbb{G}u^+|^2-| \nabla_\mathbb{G}u^-|^2=2 &\quad \mbox{ in }\
  F(u),
\\[0.3cm]u =0 &\quad \mbox{ on }\
  \partial \Omega,
\end{array}
\right.
 \end{equation}
 where $f\in L^{\infty}(\Omega)$ is a positive bounded function, we prove the existence of solutions for problem $(\ref{1})$ using a monotonicity result.\\
 The analogue of problem $(\ref{1})$ in the Euclidean case involving the Laplacian has been studied by K. Perera in \cite{perrera class}. Thus, the problem $(\ref{1})$ becomes
 \begin{equation}\label{3}
  \left\{\begin{array}{ll}  -\Delta u=\lambda \chi_{\{u>1\}}g(x,(u-1)_+) & \quad \mbox{
in }\ \Omega\setminus F(u),\\[0.3cm]| \nabla u^+|^2-| \nabla u^-|^2=2 &\quad \mbox{ in }\
  F(u),
\\[0.3cm]u =0 &\quad \mbox{ on }\
  \partial \Omega,
\end{array}
\right.
 \end{equation}
where $\Omega\subset \mathbb{R}^N, N\geq 2$ is a smooth bounded domain.\\
We should point out that many authors have investigated the special case $g(x,s)=0.$ The problem $(\ref{3})$ reduces to the following problem
\begin{equation}\label{4}
  \left\{\begin{array}{ll}  -\Delta u=0 & \quad \mbox{
in }\ \Omega\setminus F(u),\\[0.3cm]| \nabla u^+|^2-| \nabla u^-|^2=2 &\quad \mbox{ in }\
  F(u),
\\[0.3cm]u =0 &\quad \mbox{ on }\
  \partial \Omega.
\end{array}
\right.
 \end{equation}
 For more details, we refer the reader to references \cite{alt-caffarelli,alt-caffarelli-fried,caffa-annal,weiss1,weiss2}.\\
 For the case $g(x,s)=1,$ the problem $(\ref{3})$ is the well-known Prandlt-Batchelor free boundary problem arising in fluid dynamics when $N=2.$ ( See \cite{Batchelor1, Batchelor2}). Others generalizations can be founded in \cite{Caflish,Elcrat}.\\
 The study of problem $(\ref{3})$ for the suplinear case $g(x,s)=s^{p-1},$ $$\left\{\begin{array}{ll} 2<p<\infty  & \quad \mbox{if }\ N=2\\[0.3cm]
				2<p<\frac{2N}{N-2}       & \quad \mbox{if }\ N\geq 3
			\end{array}
			\right.	
		$$
related to the plasma confinement have been extensively investigated for over four decades. See for instance \cite{Frie Plasma, Fri Liu, Temam, Temam1}.\\
Note also that for a subcritical case, a mountain pass solution was obtained by Jerison and Perera in \cite{jerison} and for the critical case in \cite{yang}.\\\\
On the other hand, some few papers have considered the study of free boundary problems in non Euclidean domains. One of seminal works is Ferrari and Valdinoci \cite{valdinoc-Ferrari} where the authors prove some density estimate for local minima in the Heisenberg group. In this direction, recently, Ferrari and Forcillo \cite{Forcilli-Ferrari} prove the Alt-Caffarelli-Friedman monotonicity formula in particular form in the framework of the Heisenberg group. This fundamental tool for free boundary problems was introduced in the literature by Alt-Caffarelli-Friedman in \cite{alt-caffarelli-fried}. We cite also the interesting papers  \cite{Danielli1, Danielli2} dealing with the obstacle problem.\\\\
In this paper, motivated by many applications of the Heisenberg group in the real-world (see e.g \cite{Binz} ), we shall study a general class of free boundary problem over a stratified Lie group ( in particular the Heisenberg group). We present an extension of the result given in \cite{perrera class}. \\\\
The main result of this work is the following
 \begin{theorem}
Assume $(g_1)$ and $(g_2).$ Then, there exists a $\lambda^*>0$ such that for $0<\lambda<\lambda*,$ there exits a positive solution $u\in W_0^{1,2}(\Omega)\cap C^2(\overline{\Omega}\setminus F(u))$ satisfying the equation $ -\mathcal{L}u=\lambda \chi_{\{u>1\}}g(x,(u-1)_+)$ in the classical sense in $\Omega\setminus F(u),$ the free boundary condition in the sense of viscosity and vanishing continuously on $\partial \Omega.$
\end{theorem}
Note that the variational functional associated to problem $(\ref{1})$ is given by
$$E(u)=\int_{\Omega}\left[\frac{|\nabla_{\mathbb{G}}u|^2}{2}+\chi_{\{u>1\}}(x)-\lambda G(x,(u-1)_+)\right]dx,$$
where $G(x,s)=\int_{0}^{s}g(x,t)dt,\quad s\geq 0.$\\
We remark that $E$ is not $C^1$ du the presence of $\chi_{\{u>1\}}$ and hence the classical variational methods are not applicable. So, we use an approximation technique to resolve this issue.\\
This paper is organized as follows : In section 2, we collect some definitions and preliminary results concerning the stratified Lie group. In section 3, we recall the monotonicity lemma in the non-Euclidean setup and prove the convergence lemma. Finally, in section 4, we prove the main result of this paper.
\section{Mathematical preliminaries}
In this section, we provide the necessary tools in order to study problem $(\ref{1}).$ The reader who is familiar with these notions may go directly to the next section. We begin by giving the definitions of homogeneous Lie group.
\begin{definition}
Let $\mathbb{G}=(\mathbb{R}^N,*)$ be a Lie group on $\mathbb{R}^N.$ We say that $G$ is a homogeneous if any $\delta>0,$ there exists an automorphism $T_{\delta}:\mathbb{G}\rightarrow \mathbb{G}$ defined by
$$T_{\delta}(x)=(\delta^{r_1}x_1,\delta^{r_2}x_2,...,\delta^{r_N}x_N),\quad r_i>0,i=1,...,N\quad\hbox{for}\quad x=(x_1,...,x_N). $$
The map $T_{\delta}$ is called a dilation on $\mathbb{G}.$ The positive number $Q=\Sum_{i=1}^{N}r_i$ called homogeneous dimension of the group $\mathbb{G}.$
\end{definition}
\begin{definition}
We say that $\mathbb{G}=(\mathbb{R}^N,*)$ is a homogeneous Carnot group or a homogeneous stratified Lie group if the following properties hold\\
1) $\mathbb{R}^N$ can be split as $\mathbb{R}^N=\mathbb{R}^{N_1}\times...\mathbb{R}^{N_r}$ for some natural number $N_1,N_2,...,N_r$ such that $N_1+N_2+...+N_r=N$ and the dilation $T_{\delta}:\mathbb{R}^N\rightarrow \mathbb{R}^N,$
$$T_{\delta}(x)=(\delta^{1}x^{(1)},\delta^{2}x^{(2)},...,\delta^{r}x^{(r)}),\quad x^{(i)\in \mathbb{R}^N}$$
is an automorphism of the group $\mathbb{G}$ for every $\delta>0.$\\
2) If $N_1$ is as above, let $Z_1,...,Z_{N_1}$ be the left invariant vector fields on $\mathbb{G}$ such that $Z_i(0)=\frac{\partial}{\partial x_i}\large|_0$ for $i=1,...,N,$ then
\begin{equation}\label{5}Rank(Lie\{Z_1,...,Z_{N_1}\})(x)=N,\quad \forall x\in \mathbb{R}^N
\end{equation}
The vector fields $Z_1,...,Z_{N_1}$ are called the (Jacobian) generators of $\mathbb{G}$ whereas any basis for $span{Z_1,...,Z_{N_1}}$ is called a systm of generator of $\mathbb{G}$. We also say that $\mathbb{G}$ has step $r$ and $N_1$ generators.\\
Note also that $(\ref{5})$ is the well-known the H\"{o}rmander condition.
\end{definition}
An interesting example is the $\mathbb{H}^1=(\mathbb{R}^3,*) $ known as the Heisenberg-Weyl group on $\mathbb{R}^3$ which is a stratified Lie group of step $2$ and $2$ generators. Indeed, it is a homogeneous Lie group with dilation
$$T_{\delta}(x_1,x_2,x_3)=(\delta x_1,\delta x_2,\delta^2 x_3)$$
and for $$Z_1=\frac{\partial}{\partial x_1}+2x_2 \frac{\partial }{\partial x_3},\quad Z_2=\frac{\partial}{\partial x_2}-2 x_1\frac{\partial}{\partial x_3},$$
we have
$$rank(Lie\{Z_1,Z_2\})(x)=3,\quad \hbox{for every}\quad x\in \mathbb{R}^3.$$
Thus, $1)$ and $2)$ of the above definition are satisfied.\\\\
In the sequel, we let $N_1=N$ and let $Z=(Z_1,...,Z_N)$ be any family of vector fields on $\mathbb{R}^N.$
\begin{definition}
A Lipschitz continuous curve $\gamma:[0,T]\rightarrow \mathbb{R}^N, T\geq 0$ is $Z-$admissible if there exists a bounded functions $h=(h_1,h_2,...,h_N):[0,T]\rightarrow\mathbb{R}^N$ such that
$$\gamma'(t)=\Sum_{j=1}^N h_j Z_j(\gamma(t)),\quad \hbox{for a.e}\quad t\in [0,T].$$
In general, the function $h$ is not unique.
\end{definition}
Then, the following definition makes sense.
\begin{definition}(Carnot-Carathéodory distance)
For any $x,y\in \mathbb{G},$ the Carnot-Carathéodory distance is defined by
$$d_Z(x,y)=inf\{T\geq 0, \hbox{then there exits an admissible}\quad \gamma:[0,T]\rightarrow \mathbb{G}\ \hbox{such that}\quad \gamma(0)=x,\gamma(T)=y\}.$$
If the above set is empty put $d_Z(x,y)=0.$\\
If $d_Z(x,y)<+\infty$ for all $x,y\in \mathbb{G},$ then $(\mathbb{G},d_Z))$ is a metric space known as the Carnot-Carathéodory space.
\end{definition}
Now, we define the Sub-Laplacian ( or horizontal Laplacian) on $\mathbb{G}$ as
$$\mathcal{L}:=Z_1^2+...+Z_N^2.$$
The horizontal gradient on $\mathbb{G}$ is defined as
$$\nabla_\mathbb{G}:=(Z_1,...,Z_N),$$
and the horizontal divergence on $\mathbb{G}$ is defined by
$$div_\mathbb{G}v=\nabla_\mathbb{G}.v.$$
Thus, the sub-Laplacian $\mathcal{L}$ can be written as
$$\mathcal{L}u=div_\mathbb{G}(\nabla_\mathbb{G}u).$$
On the other hand, for $x\in \mathbb{G},$ the Haar measure ( denoted by $dx$) is the Lebesgue measure on $\mathbb{R}^N.$ So, for $1\leq p\leq \infty,$ the space $L^p(\mathbb{G})$ or simply $L^p$ denote the usual Lebesgue space on $\mathbb{G}$ with respect to the Haar measure with the  norm
$$||.||_{L^p}=\|.\|_{L^p(\mathbb{G})}=\|.\|_p$$
given  for $p\in [1,+\infty)$ by
$$\|f\|_p=\left(\int_\mathbb{G}|f(x)|^p dx\right)^{1/p},$$
and for $p=\infty$ by
$$\|f\|_{\infty}=\Sup_{x\in \mathbb{G}}|f(x)|.$$
Here, the supremum with respect to the Haar measure. We recall that the Haar measure for a measurable subset of $\mathbb{G}$ is denoted by $\mu(.)$ and satisfies the properties
$$\mu(T_{\delta}(M))=\delta^Q\mu(M)$$
where $M\subset \mathbb{G}$ is a measurable set.\\\\
For the variational setting of our problem, we define the Sobolev space $W^{1,p}(\Omega)$ where $\Omega$ is an open subset on a stratified Lie group as
$$W^{1,p}(\Omega)=\{u\in L^p(\Omega):\quad |\nabla_\mathbb{G}u|\in L^p(\Omega)\}.$$
A norm on this space
$$\|u\|_{1,p}=\|u\|_p+\|\nabla_\mathbb{G} u\|_{p}.$$
However, the Sobolev space $W_0^{1,p}(\Omega)$ is defined as
$$W_0^{1,p}(\Omega)=\{u\in W^{1,p}(\Omega),\quad u=0\quad \hbox{on}\quad \partial \Omega\},$$
where $u=0$ on $\partial \Omega$ is the usual trace sense. This space is separable and uniformly convex Banach space.\\
If $\Omega$ is bounded, then for $1\leq p\leq Q,$ $W_0^{1,p}(\Omega)$ is continuously embedded in $L^q$ for $1\leq q\leq p^*=\frac{Qp}{Q-p}$ and compact if $1\leq q< p^*.$\\
Also, the Sobolev space $W_{loc}^{1,p}(\Omega)$ is the space of all functions in $W^{1,p}(K)$ for any compact set $K\subset \Omega.$\\
Hence, throughout  this work, we need the following result
\begin{lemma}\cite{Repsov}
  For $u\in W_{loc}^{1,p}(\Omega),$ the Radon measure $\mathcal{M}=\mathcal{L}u$ is nonnegative and supported on $\Omega\cap\{u<1\}.$
\end{lemma}
We refer the reader to \cite{D1,D2} for more details.\\\\
Finally, we state the analogue of divergence theorem in the Euclidean setup for sub-Laplacian on stratified Lie group. This result is very important for the proofs of our results.
\begin{proposition}\cite{Ruzhansky}
Let $f_n\in C^1(\Omega)\cap C(\overline{\Omega}),$ $n=1,2,...,N.$ Then, for each $n,$ we have
$$\int_{\Omega}Z_n f_n d\nu=\int_{\partial \Omega}f_n<Z_n,d\nu>.$$
Consequently,
$$\int_{\Omega}\Sum_{n=1}^{N}Z_n f_n d\nu=\int_{\partial \Omega} \Sum_{n=1}^{N} f_n<Z_n,d\nu>.$$
\end{proposition}
\section{Approximate solution}
In order to give the proof of the main result, we recall the important monotonicity lemma in a non-Euclidean setup proved by Choudhuri and Repov$\breve{s}$ in \cite{Repsov}. The classical version ( in the Euclidean setting) have been proved by Caffarelli et al in \cite{caffa-annal}.
\begin{lemma}\cite{Repsov}
Let $u>0$ be a Lipschitz continuous function on the unit ball $B_1(0)\subset \mathbb{G}$ satisfying the distributional inequalities
$$\pm \mathcal{L}u\leq \left(\frac{\lambda}{\varepsilon}\chi_{\{|u-1|<\varepsilon\}}(x)\mathfrak{F}(|\nabla_\mathbb{G}u|)+A\right),$$
for constants $A>0,$ $0<\varepsilon<1.$ Suppose further that $\mathfrak{F}$ is a continuous function such that $\mathfrak{F}(t)=o(t^2)$ near infinity. Then, there exist $C=C(N,A)>0$ and $\int_{B_1(0)}u^2 dx$ but not on $\varepsilon$ such that
$$\esssup_{x\in B_{1/2}(0)}\{|\nabla_\mathbb{G}u|\}\leq C.$$
\end{lemma}

Now,  we approximate $E$ by $C^1-$functionals. So, let $\beta:\mathbb{R}\rightarrow [0,2]$ be a function such that\\
$(1)$ $\beta\in C^{\infty}(\mathbb{R}),$ support of $\beta$ lies in $[0,1]$ and it is positive in $(0,1).$\\
$(2)$ $\int_{0}^{1}\beta(s)ds=1.$\\\\
Then, let $B(s)=\int_{0}^{s}\beta(t)dt$ where $B:\mathbb{R}\rightarrow [0,1]$ is a smooth and nondecreasing function such that
$$B(s)=0, \quad \hbox{for}\quad s\leq 0,\quad 0<B(s)<1\quad \hbox{for}\quad 0<s<1\quad \hbox{and}\quad B(s)=1\quad \hbox{for}\quad s\geq 1.$$
Hence, for $\varepsilon>0,$ let
$$E_{\varepsilon}(u)=\int_{\Omega}\left[\frac{|\nabla_\mathbb{G}|^2}{2}+B(\frac{u-1}{\varepsilon})-\lambda G_{\varepsilon}(x,(u-1)_+)\right]dx$$
where $G_{\varepsilon}(x,s)=\int_{0}^{s}g_{\varepsilon}(x,t)dt,\quad s\geq 0,$ \quad $g_{\varepsilon}(x,s)=B(\frac{s}{\varepsilon})g(x,s).$\\\\
The function $E_{\varepsilon}$ is of class $C^1$ and its critical points are the solutions of the following problem
\begin{equation}\label{app}
  \left\{\begin{array}{ll}  -\mathcal{L}u=-\frac{1}{\varepsilon}\beta(\frac{u-1}{\varepsilon})+\lambda g_{\varepsilon}(x,(u-1)_+) & \quad \mbox{
in }\ \Omega,
\\[0.3cm]u =0 &\quad \mbox{ on }\
  \partial \Omega.
\end{array}
\right.
 \end{equation}
 Let $\varepsilon_j\searrow 0,$ let $u_j$ be a critical point of $E_{\varepsilon_j}.$ We prove the following convergence result
 \begin{lemma}
 Assume $g_1)$ and $g_2).$ If $(u_j)$ is bounded in $W_0^{1,2}(\Omega)\cap L^{\infty}(\Omega),$ then there exists a Lipschitz continuous function $u$ on $\overline{\Omega}$ such that
 \begin{enumerate}[label=(\roman*)]
        \item $u_j\rightarrow u$ uniformly over $\overline{\Omega}.$
        \item $u_j\rightarrow u$ locally in $C^1(\overline{\Omega}\setminus \{u=1\}).$
        \item $u_j\rightarrow u$ strongly in $W_0^{1,2}(\Omega).$
        \item $E(u)\leq \Lim \inf E_{\varepsilon_j}(u_j)\leq \Lim \sup E_{\varepsilon_j}(u_j)\leq E(u)+\left|\{u=1\}\right|.$
    \end{enumerate}
 Moreover, $u$ satisfies $-\mathcal{L}u=\lambda \chi_{\{u>1\}}g(x,(u-1)_+)$ in the classical sense in $\Omega\setminus F(u)$ and $u$ satisfies the free boundary condition in the generalized sense and vanishes continuously on $\Omega.$
 \end{lemma}
 {\bf Proof of Lemma 3.1}
 Assume that $0<\varepsilon_j<1.$ So, $u_j$ is a solution of
 \begin{equation}\label{appapp}
  \left\{\begin{array}{ll}  -\mathcal{L}u_j=-\frac{1}{\varepsilon}\beta(\frac{u_j-1}{\varepsilon})+\lambda g_{\varepsilon}(x,(u_j-1)_+) & \quad \mbox{
in }\ \Omega,
\\[0.3cm]u_j =0 &\quad \mbox{ on }\
  \partial \Omega.
\end{array}
\right.
 \end{equation}
 Since $(u_j)$ is bounded in $L^{\infty}(\Omega),$ then $0\leq g_{\varepsilon_j}(x,(u_j-1)_+)\leq A_0$ for some constant $A_0>0.$\\\\
 Let $\varphi_0$ be a solution of
 \begin{equation}\label{pb A}
  \left\{\begin{array}{ll}  -\mathcal{L}\varphi_0=\lambda A_0 & \quad \mbox{
in }\ \Omega,
\\[0.3cm]\varphi_0 =0 &\quad \mbox{ on }\
  \partial \Omega.
\end{array}
\right.
 \end{equation}
 Since $\beta \geq 0,$ then $-\mathcal{L}u_j\leq \lambda A_0.$ So, by the maximum principle, we have
 $$0\leq u_j\leq \varphi_0(x),\quad \quad \forall x\in \Omega.$$
 Thus, taking into account that $\{u_j\geq 1\}\subset \{\varphi_0\geq 1\},$ it follows that $\varphi_0$ gives a uniform lower bound, say $d_0$ on the distance from the set $\{u_j\geq 1\}$ to $\partial \Omega.$\\\\
 Also, since the sequence $(u_j)$ is bounded with respect to $C^{2,\alpha}-$norm, so, from standard boundary regularity theory, it has  a convergent subsequence in $C^2-$norm on $\frac{d_0}{2}$ neighborhood of $\partial \Omega.$\\\\
 Observe that $0\leq \beta\leq 2\chi_{(-1,1)}, $ then
 $$\pm \mathcal{L}u_j=\pm \frac{1}{\varepsilon_j}\beta(\frac{(u_j-1)_+}{\varepsilon_j})\pm \lambda g_{\varepsilon_j}(x,(u_j-1)_+)$$
 $$\leq \frac{2}{\varepsilon_j}\chi_{\{|u_j-1|<\varepsilon_j\}}(x)+\lambda A_0.$$
 Using Lemma $3.1,$ then there exists $A>0$ such that
 $$\esssup_{x\in B_{r/2}(x_0)}\{|\nabla_\mathbb{G}u_j|\}\leq \frac{A}{r},\quad \hbox{for}\quad r>0,$$
 for which $B_r(x_0)\subset \Omega.$\\
 Because $(u_j)$ is a sequence of Lipschitz continuous functions, then we have
 $$\Sup_{x\in B_{r/2}(x_0)}\{|\nabla_\mathbb{G}u_j|\}\leq \frac{A}{r}.$$
 Therefore, the family $(u_j)$ is uniformly Lipschitz continuous on the compact subsets, say $K$ such that $d(K,\partial \Omega)\geq \frac{d_0}{2}.$\\
 The application of Ascoli-Arzela gives a subsequence namely also $(u_j)$ that converge uniformly to a Lipschitz continuous function $u$ in $\Omega$ and finally, the Banach-Alaoglu, we have $u_j\rightharpoonup u$ in $W_0^{1,2}.$\\
 Now, we show that $u$ satisfies
 \begin{equation}\label{et}
 -\mathcal{L}u=\lambda \chi_{\{u>1\}}g(x,(u-1)_+)\quad \hbox{in}\quad \{u\neq 1\}.
 \end{equation}
 Let $\varphi\in C_0^{\infty}(\{u>1\}).$ Then, on the support of $\varphi,$ $u\geq 1+2\delta$ for $\delta>0.$\\
For large $j,$ $\delta_j<\delta$ and  $|u_j-u|<\delta,$ we have $u_j\geq 1+\delta_j$ on the support of $\varphi.$\\
So, testing $(\ref{appapp})$ with $\varphi$ gives
$$\int_{\Omega}\widetilde{\nabla}u_j\varphi dx=\lambda\int_{\Omega}g(x,u_j-1)\varphi dx,$$
where $\widetilde{\nabla}\eta=\Sum_{k=1}^{N}X_k\eta X_k.$\\
Passing to the limit $j\rightarrow +\infty,$ we have
$$\int_{\Omega}\widetilde{\nabla}u\varphi dx=\lambda\int_{\Omega}g(x,u-1)\varphi dx.$$
Hence $u$ is a distributional (and thus a classical) solution of
\begin{equation}\label{E1}
 -\mathcal{L}u=\lambda g(x,(u-1)_+)\quad \hbox{in}\quad \{u> 1\}.
\end{equation}
By the same argument, taking $\varphi\in C_0^{\infty}(\{u<1\}),$ we find $\delta>0$ such that $u\leq 1-2\delta$ giving also $u_j<1-\delta.$\\
testing $(\ref{appapp})$ with any nonnegative $\varphi$ and passing to the limit yields
$$-\mathcal{L}u\leq \lambda g(x,(u-1)_+)\quad \hbox{in}\quad \Omega$$
in the distributional sense. So,
$\mathcal{L}(u-1)_{-}=\mu$ is a positive Radon measure supported on $\Omega \cap \{u<1\}$ (see Lemma $2.1$ in section 2).\\
By the classical result of regularity ( see section 9.4 in \cite{gilbarg}), we obtain that $u\in W_{loc}^{2,2}(\{u\leq 1\}^o)$ and hence $\mu$ is actually supported on $\Omega\cap \partial \{u<1\}\cap \{u>1\},$ So $u$ satisfies $\mathcal{L} u=0$ on the set $\{u\leq 1\}^o.$\\\\
Now, we prove ii). We have proved that $u_j\rightarrow u$ with respect to $u$ with respect to $C^2$ norm in a neighborhood of $\partial \Omega$ in $\overline{\Omega}.$ Then, for $\delta>0,$ we have $u\geq 1+2\delta$ in the set $U\subset\subset\{u>1\}.$\\
So, for large $j$ with $\delta_j<\delta,$ we have $|u_j-u|<\delta$ in $\Omega$ and hence $u_j\geq 1+\delta_j$ in $U.$\\
From $(\ref{app}),$ we have
$$-\mathcal{L} u_j=\lambda g(x,(u_j-1))\quad \hbox{in}\quad U.$$
Observe that $g$ is locally H\"{o}lder continuous and $u_j\rightarrow u$ uniformly, then $g(x,(u_j-1))\rightarrow g(x,(u-1))$ in $L^p(U),$ $p>1.$ Also, since $W^{2,2}(U)\hookrightarrow C^1(U),$ then, we have $u_j\rightarrow u$ in $C^1(U)$ which implies that $u_j\rightarrow u$ in $C^1(\{u>1\}).$\\
The same method can be applied to prove that $u_j\rightarrow u$ in $C^1(\{u<1\}).$\\\\
To prove iii), we remark that $u_j\rightharpoonup u$ in $W_0^{1,2}(\Omega)$ and by the weak lower semicontinuity of the norm $\|.\|,$ we have
$$\|u\|\leq \varliminf \|u_j\|.$$
Hence, to show iii), it suffices  to prove that $\varlimsup \|u_j\|\leq \|u\|.$\\
So, we multiply $(\ref{app})$ by $(u_j-1)$ and then integrate by parts. Using also the fact that $\tan(\frac{t}{\delta_j})\geq 0$ for all $t>0.$ Then, we obtain
\begin{align}\int_{\Omega}|\nabla_\mathbb{G}|^2dx & \leq \lambda \int_{\Omega}g(x,(u_j-1)_+)(u_j-1)_+dx-\int_{\partial \Omega}u_j<X_i,dn>ds \nonumber\\
&\rightarrow \lambda \int_{\Omega}g(x,(u-1)_+)(u-1)_+dx-\int_{\partial \Omega}u<X_i,dn>ds
\end{align}
where $n$ is the outward unit normal to $\partial \Omega.$\\
Now, fix $0<\varepsilon<1$ and testing the equation $(\ref{E1})$ with $\varphi=(u-1-\varepsilon)_+$ yields
\begin{equation}\label{I1}
\int_{\{u>1+\varepsilon\}}|\nabla_\mathbb{G}u|^2dx=\int_{\Omega}\lambda g(x,(u-1)_+)(u-1-\varepsilon)_+dx.
\end{equation}
On the other hand, integrating $(u-1+\varepsilon)\Delta u=0$ over $\Omega$ gives
\begin{equation}\label{I2}
\int_{\{u<1-\varepsilon\}}|\nabla_\mathbb{G}u|^2dx=-(1-\varepsilon)\int_{\partial\Omega}u<X_i,dn>ds.
\end{equation}
Adding, $(\ref{I1})$ and $(\ref{I2})$ and letting $\varepsilon\rightarrow 0$ gives
$$\int_{\Omega}|\nabla_\mathbb{G}u|^2dx=\int_{\Omega}\lambda g(x,(u-1)_+)(u-1)_+dx-\int_{\partial \Omega}u<X_i,dn>ds$$
since $\int_{\{u=1\}}|\nabla_\mathbb{G}u|^2dx=0.$ So, using $(11),$ we have
$$\varlimsup \int_{\Omega}|\nabla_\mathbb{G}u_j|^2dx\leq \int_{\Omega}|\nabla_\mathbb{G}u|^2dx$$
as wanted.\\
For the last statement $iv),$ we have
$$E_{\varepsilon_j}(u_j)=\int_{\Omega}\left[\frac{1}{2}|\nabla_\mathbb{G}u_j|^2+B\left(\frac{u_j-1}{\varepsilon_j}\right)\chi_{\{u\neq 1\}}-\lambda G_{\varepsilon_j}(x,(u_j-1)_+)\right]dx+\int_{\{u=1\}}B\left(\frac{u_j-1}{\varepsilon_j}\right)dx.$$
Since $B\left(\frac{u_j-1}{\varepsilon_j}\right)\chi_{\{u\neq 1\}}$ converges pointwise to $\chi_{\{u>1\}}$ and is bounded by $1$
 and $G_{\varepsilon_j}(x,(u_j-1)_+)$ converges to $G(x,(u-1)_+),$ then the first integral converges to $E(u).$\\
 On the other hand, observe that
 $$0\leq \int_{\{u=1\}}B\left(\frac{u_j-1}{\varepsilon_j}\right)dx\leq \left|\{u=1\}\right|. $$
 Thus, $iv)$ follows.\\\\
 Now, to conclude the proof of Lemma $3.1,$ we choose $\overrightarrow{\varphi}\in C_0^1(\Omega,G)$ such that $u\neq 1$ a.e on the support of $\overrightarrow{\varphi}.$ \\
 Multiplying by $\Sum_{k=1}^{N}\varphi_k X_k u_n$ the weak formulation of $(\ref{appapp})$ and integrate over the set $\{1-\varepsilon^-< u_n<1+\varepsilon^+\},$ we get
 \begin{align}\label{I3}
   \int_{\{1-\varepsilon^-< u_n<1+\varepsilon^+\}}\left[-\Delta_\mathbb{G}u_n+\frac{1}{\varepsilon_n}\beta(\frac{u_n-1}{\varepsilon_n})\right]\Sum_{k=1}^{N}\varphi_k X_k u_n dx & = \nonumber\\
     =\lambda\int_{\{1-\varepsilon^-< u_n<1+\varepsilon^+\}}g_{\varepsilon_n}(x,(u_n-1)_+)\Sum_{k=1}^{N}\varphi_k X_k u_ndx.
 \end{align}
 The term of the left in $(\ref{I3})$ can be simplified as
\begin{align}
\nabla_\mathbb{G}\left(\frac{1}{2}|\nabla_\mathbb{G}u_n|^2\overrightarrow{\varphi}-\left(\Sum_{k=1}^{N}X_k u_n \varphi_k\right)|\nabla_\mathbb{G}u_n|\right)+\Sum_{k=1}^{N}\Sum_{l=1}^{N}X_l \varphi_k u_n X_k u_n \nonumber\\-\frac{1}{2}|\nabla_\mathbb{G}u_n|^2\nabla_\mathbb{G}\overrightarrow{\varphi}
+\Sum_{k=1}^{N}\varphi_k X_k B(\frac{u_n-1}{\varepsilon_n}).
\end{align}
We integer by part to obtain
\begin{align}\label{I4}
  \int_{\{1-\varepsilon^-< u_n<1+\varepsilon^+\}}\frac{1}{2}|\nabla_\mathbb{G}u_n|^2\Sum_{k=1}^{N}\varphi_k<X_k,dn>-(\Sum_{k=1}^{N}X_k u_n \varphi_k)\Sum_{l=1}^{N}X_l u_n<X_l,dn> \nonumber  \\
  + \int_{\{1-\varepsilon^-< u_n<1+\varepsilon^+\}}B(\frac{u_n-1}{\varepsilon_n})\Sum_{k=1}^{N}\varphi_k<X_k,dn> \nonumber  \\
 = \int_{\{1-\varepsilon^-< u_n<1+\varepsilon^+\}}(\frac{1}{2}|\nabla_\mathbb{G}|^2\Sum_{k=1}^{N} X_k\varphi_k-\Sum_{k=1}^{N}\Sum_{l=1}^{N}X_k\varphi_lX_l u_n X_k u_n)\nonumber \\
 + \int_{\{1-\varepsilon^-< u_n<1+\varepsilon^+\}}\left[B(\frac{u_n-1}{\varepsilon_n})\Sum_{k=1}^{N}X_k \varphi_k+\lambda g_{\varepsilon_n}(x,(u_n-1)_+)\Sum_{k=1}^{N}X_k \varphi_k\right]dx.
\end{align}
The integral on the left converges to
\begin{align}
& \int_{\left\{u=1+\epsilon^{+}\right\} \cup\left\{u=1-\epsilon^{-}\right\}}\left[\frac{1}{2}\left|\nabla_{\mathbb{G}} u\right|^{2} \sum_{k=1}^{N} \varphi_{k}\left\langle X_{k}, d n\right\rangle-\left(\sum_{k=1}^{N} X_{k} u \varphi_{k}\right) \sum_{l=1}^{N} X_{l} u\left\langle X_{l}, d n\right\rangle\right. \\
+ & \left.\int_{\left\{u=1+\epsilon^{+}\right\}} \sum_{k=1}^{N} \varphi_{k}\left\langle X_{k}, d n\right\rangle\right] \nonumber \\
= & \int_{\left\{u=1+\epsilon^{+}\right\} \cup\left\{u=1-\epsilon^{-}\right\}}\left[\left(1-\frac{1}{2}\left|\nabla_{\mathbb{G}} u\right|^{2}\right) \sum_{k=1}^{N} \varphi_{k}\left\langle X_{k}, d n\right\rangle-\sum_{k \neq l ; 1 \leq k, l \leq N} \varphi_{k} X_{l} u X_{k} u\left\langle X_{l}, d n\right\rangle\right]  \\
= & \int_{\left\{u=1+\epsilon^{+}\right\}}\left[\left(1-\frac{1}{2}\left|\nabla_{\mathbb{G}} u\right|^{2}\right) \sum_{k=1}^{N} \varphi_{k}\left\langle X_{k}, d n\right\rangle-\sum_{k \neq l ; 1 \leq k, l \leq N} \varphi_{k} X_{l} u X_{k} u\left\langle X_{l}, d n\right\rangle\right]\nonumber
\end{align}
$$
\begin{aligned}
& -\int_{\left\{u=1-\epsilon^{-}\right\}}\left[\left(\frac{1}{2}\left|\nabla_{\mathbb{G}} u\right|^{2}\right) \sum_{k=1}^{N} \varphi_{k}\left\langle X_{k}, d n\right\rangle-\sum_{k \neq l ; 1 \leq k, l \leq N} \varphi_{k} X_{l} u X_{k} u\left\langle X_{l}, d n\right\rangle\right] \\
= & \int_{\left\{1-\epsilon^{-}<u<1+\epsilon^{+}\right\}}\left(\frac{1}{2}\left|\nabla_{\mathbb{G}} u\right|^{2} \sum_{k=1}^{N} X_{k} \varphi_{k}-\sum_{k=1}^{N} \sum_{l=1}^{N} X_{k} \varphi_{l} X_{l} u X_{k} u\right) d x \\
& \int_{\left\{1-\epsilon^{-}<u<1+\epsilon^{+}\right\}}\left[\sum_{k=1}^{N} X_{k} \varphi_{k}+\lambda g(x,(u-1)_+) f \sum_{k=1}^{N} X_{k} \varphi_{k}\right] d x,
\end{aligned}
$$

as $n \rightarrow \infty$.\\
Observe that the normal vector at the point $P$ on the set $\left\{u=1+\epsilon^{+}\right\} \cup\left\{u=1-\epsilon^{-}\right\}$is $n= \pm \frac{\nabla_{\mathbb{G}} u(P)}{\left|\nabla_{\mathbb{G}} u(P)\right|}$. So, equation $(18)$ becomes when $\varepsilon\rightarrow 0$
$$0=\lim _{\epsilon \rightarrow 0} \int_{\left\{u=1+\epsilon^{+}\right\}}\left[\left(1-\frac{1}{2}\left|\nabla_{\mathbb{G}} u\right|^{2}\right) \sum_{k=1}^{N} \varphi_{k}\left\langle X_{k}, d n\right\rangle\right] \\
-\lim _{\epsilon \rightarrow 0} \int_{\left\{u=1-\epsilon^{-}\right\}}\left[\left(\frac{1}{2}\left|\nabla_{\mathbb{G}} u\right|^{2}\right) \sum_{k=1}^{N} \varphi_{k}\left\langle X_{k}, d n\right\rangle\right] .$$
Hence, we conclude that $u$ satisfies the free boundary condition in the sense of viscosity.
 \section{Proof of the main result}
 In this section, we prove our main result given in theorem 1.1. First, taking account that $g$ verifies $g_1),$ we show that $E_{\varepsilon}$ satisfies the Palais-Smale (PS) condition. We have
 \begin{lemma}
 Assume $g_1).$ The functional $E_{\varepsilon}$ satisfies the (PS) condition.
 \end{lemma}
 {\bf Proof of Lemma 4.1.} We note that  by $g_1),$ we have
 $$E_{\varepsilon}(u_n)\geq \int_{\Omega}\left[\frac{1}{2}|\nabla_\mathbb{G} u_n|^2-\lambda(a_0(u_n-1)_++\frac{a_1
 }{p}(u_n-1)_+^p)\right]dx.$$
 If we denote by $||u||:=\int_{\Omega}|\nabla_\mathbb{G} u|dx,$ then
$$<E_{\varepsilon}'(u_n),u_n>\leq ||u_n||^2-\lambda \int_{\Omega}g(x,(u_n-1)_+)dx+\frac{2}{\varepsilon}|\Omega|. $$
A standard argument implies that $(u_n)$ is bounded in $W_{0}^{1,2}(\Omega).$ Hence, $u_n\rightharpoonup u$ in $W_{0}^{1,2}(\Omega).$\\
On the other hand, since $<E_{\varepsilon}'(u_n),u_n>\rightarrow 0$ as $n\rightarrow +\infty,$ we have
$$\Lim_{n\rightarrow +\infty}\int_{\Omega} \widetilde{\nabla} u_n v=\Lim_{n\rightarrow +\infty}\left[\int_{\Omega}\frac{1}{\varepsilon}\beta(\frac{u_n-1}{\varepsilon})vdx+\lambda\int_{\Omega}g(d,(u_n-1)_+)vdx\right],\quad \forall v\in W_{0}^{1,2}(\Omega).$$
Taking, $v=u_n-u$,  we obtain
$$\Lim_{n\rightarrow +\infty}\int_{\Omega} \widetilde{\nabla} u_n v=\Lim_{n\rightarrow +\infty}\left[\int_{\Omega}\frac{1}{\varepsilon}\beta(\frac{u_n-1}{\varepsilon})(u_n-u)dx+\lambda\int_{\Omega}g(d,(u_n-1)_+)(u_n-u)dx\right]=0.$$
So, $u_n\rightarrow u$ in $W_{0}^{1,2}(\Omega)$ and $E_{\varepsilon}$ satisfies the (PS)
 condition.$\blacksquare$\\\\
 {\bf Proof of theorem 1.1.}\\
 First, we remark that $E$ is bounded from below.\\
 By the assumption $g_2),$ $G(x,s)>0,$ $\forall x\in \Omega,$ $\forall s>0.$ Hence, for $u\in W_{0}^{1,2}(\Omega)$ with $u>1$ on a set of positive measure
 $$\int_{\Omega}G(x,(u-1)_+)dx>0.$$
 So, $E(u)\rightarrow -\infty$ as $\lambda\rightarrow \infty.$\\
 Furthermore, there exists a function $u_0\in W_{0}^{1,2}(\Omega)$ such that
 $$E_{\varepsilon}(u_0)<0=E_{\varepsilon}(0).$$
 Therefore, the class of paths
 $$\Gamma_{\varepsilon}=\{\gamma\in C([0,1],W_{0}^{1,2}(\Omega)),\gamma(0)=0, E_{\varepsilon}(\gamma(1))<0\}$$ is nonempty.\\
 By the mountain pass theorem, we have
 $$c_{\varepsilon}:=\Inf_{\gamma\in \Gamma_{\varepsilon}}\Max_{u\in \gamma([0,1])}E_{\varepsilon}(u).$$
 Since $B(\frac{t-1}{\varepsilon})\leq \chi_{\{t>1\}},$ $\forall t,$
 $$E_{\varepsilon}(u)\leq E(u),\quad \forall u\in W_{0}^{1,2}(\Omega).$$
 So, $\Gamma\subset\Gamma_{\varepsilon}$ and
 $$c_{\varepsilon}\leq \Max_{u\in \gamma([0,1])}E_{\varepsilon}(u)\leq \Max_{u\in \gamma([0,1])}E(u).$$
 Thus, $$c_{\varepsilon}\leq c:=\Inf_{\gamma\in \Gamma}\Max_{u\in \gamma([0,1])}E(u).$$
 Since $E_{\varepsilon}$ satisfies the (PS)
 condition, it follows that $E$ has a minimizer $u_{\varepsilon}$ satisfying $E_{\varepsilon}(u_{\varepsilon})=c_{\varepsilon}.$ \\
 From lemma 3.2, we know that $(u_{\varepsilon_n})$ denoted by $(u_n)$ converge uniformly in $\overline{\Omega},$ locally in $C^1(\overline{\Omega}\subset \{u=1\})$ and satisfy in $W_{0}^{1,2}(\Omega)$ to a lcally Lipschitz function $u\in W_{0}^{1,2}(\Omega)\cap C^2(\overline{\Omega}\setminus F(u)).$\\
 So, $u$ verifies the equation $-\mathcal{L} u=\lambda \chi_{\{u>1\}}g(x,(u-1)_+)$ in the classical sense in $\Omega\setminus F(u)$ and the free boundary condition in the sense of viscosity and vanishing continuously on $\partial \Omega.$
\section*{Ethical Approval}
This declaration is not applicable.
\section*{Competing interests}
The authors declare that there are no competing interests between them.
\section*{ Authors' contributions }
The authors declare that they read and approved the final manuscript.
\section*{ Funding}
The authors received no direct funding for this work.
\section*{ Availability of data and materials}
This declaration is not applicable.


\begin{thebibliography}{99}
\bibitem{alt-caffarelli}H.W. Alt,  L.A. Caffarelli,  Existence and regularity for a minimum problem
with free boundary, {\em J. Reine Angew. Math,} 325, (1981), 105–144.

\bibitem{alt-caffarelli-fried}H.W. Alt, L.A. Caffarelli, A. Friedman, Variational problems with two phases and their free boundaries, {\em Trans. Am. Math. Soc}, 282(2), (1984), 431–461.

\bibitem{Batchelor1} G.K. Batchelor,  On steady laminar flow with closed streamlines at large
Reynolds number, {\em J. Fluid Mech,} 1, (1956), 177–190.

\bibitem{Batchelor2}G.K. Batchelor,  A proposal concerning laminar wakes behind bluff bodies at
large Reynolds number, {\em J. Fluid Mech,} 1, (1956), 388–398.

\bibitem{Binz}E. Binz, S. Pods, The Geometry of Heisenberg Groups: With Applications in Signal Theory, Optics, Quantization, and Field Quantization, with an appendix by Serge Preston, Math. Surveys Monogr., vol. 151,  American Mathematical Society, Providence, RI, 2008, xvi+299 pp.

    \bibitem{Frie Plasma}L. A. Caffarelli, A. Friedman,  Asymptotic estimates for the plasma problem, {\em Duke Math. J,} 47(3),
705–742, (1980).

\bibitem{caffa-annal}L.A. Caffarelli, D. Jerison, C.E. Kenig, Some new monotonicity theorems with applications to free boundary problems, {\em Ann. Math,}  155(2), (2002), 369–404.

    \bibitem{caff-salsa} L. Caffarelli, S. Salsa. A geometric approach to free boundary prob-
lems, volume 68 of Graduate Studies in Mathematics. American Mathematical
Society, Providence, RI, 2005.

\bibitem{Caflish} R. E. Caflisch,  Mathematical analysis of vortex dynamics. In: Mathematical Aspects
of Vortex Dynamics (Leesburg, VA, 1988), pp. 1–24. SIAM, Philadelphia,
PA (1989)

\bibitem{Repsov}D. Choudhuria, D.D. Repovš, On semilinear equations with free boundary conditions on stratified Lie groups, {\em J. Math.Anal.Appl,} 518, (2023), 126677.
    \bibitem{Repsov1}D. Choudhuri, L. S. Tavares and J. A.  López, A study of a critical hypoelliptic problem in a stratified Lie group, {\em Complex Variables and
Elliptic Equations,} DOI: 10.1080/17476933.2024.2310217

\bibitem{D1}    D. Danielli, A compact embedding theoremfor a class of degenerate Sobolev spaces,{\em Rend Sem
Mat Univ Poi Torino,} 49(3), 1991, 399–420.
\bibitem{D2} D. Danielli, Regularity at the boundary for solutions of nonlinear subelliptic equations,{\em Indiana
Univ Math J,} , 44(1), 1995, 269–286. Doi: 10.1512/iumj.1995.44.1988

\bibitem{Danielli1}D. Danielli, N. Garofalo, A. Petrosyan, The sub-elliptic obstacle problem: C1,a regularity of the free boundary in
Carnot groups of step two, {\em Adv. Math,} 211, no. 2, (2007),  485–516.

\bibitem{Danielli2}D. Danielli, N. Garofalo, S. Salsa, Variational inequalities with lack of ellipticity. I. Optimal interior regularity and
non-degeneracy of the free boundary, {\em Indiana Univ. Math. J,} 52, no. 2, (2003),  361–398.

\bibitem{Elcrat}A.R. Elcrat, K. G.  Miller,  Variational formulas on Lipschitz domains, {\em Trans.
Am. Math. Soc,} 347(7), 2669–2678, (1995).

\bibitem{Forcilli-Ferrari}F. Ferrari, N. Forcillo, A new glance to the Alt-Caffarelli-Friedman monotonicity formula, {\em Mathematics in Engineering,} 2(4), 2020, 657-679. Doi: 10.3934/mine.2020030

\bibitem{valdinoc-Ferrari}F. Ferrari, E. Valdinoci, Density estimates for a fluid jet model in the Heisenberg group, {\em J. Math. Anal. Appl,} 382(1) (2011) 448–468.
  \bibitem{Frie-Livre}A. Friedman. Variational principles and free-boundary problems. Robert E.
Krieger Publishing Co. Inc., Malabar, FL, second edition, 1988.

\bibitem{Fri Liu} A. Friedman, Y. Liu,  A free boundary problem arising inmagnetohydrodynamic system, {\em Ann. Scuola
Norm. Sup. Pisa Cl. Sci,} (4) 22(3), 375–448, (1995).

 \bibitem{gilbarg}  D. Gilbarg, N.S. Trudinger, Elliptic Partial Differential Equations of Second Order, Springer-Verlag, Berlin, Heidelberg, 2001.
\bibitem{jerison}D. Jerison, K. Perera, Higher critical points in an elliptic free boundary problem, {\em J. Geom. Anal,} 28(2) (2018) 1258–1294.

\bibitem{KINDER} D. Kinderlehrer, L. Nirenberg, Regularity in free boundary problems, {\em Ann.
Scuola Norm. Sup. Pisa Cl. Sci,} (4), 4(2):373–391, 1977.
\bibitem{perrera class}K. Perera, On a class of elliptic free boundary problems with multiple solutions, {\em Nonlinear Differ. Equ. Appl,} 28, (2021), 36.
\bibitem{Ruzhansky}M. Ruzhansky, D. Suragan, Layer potentials, Kac’s problem, and refined Hardy inequality on homogeneous Carnot groups, {\em Adv. Math,} 308, (2017), 483–528.

 \bibitem{Temam}   R. Temam,  A non-linear eigenvalue problem: the shape at equilibrium of a confined plasma,{\em Arch.
Rational Mech. Anal,} 60(1), 51–73, (1975/76).

\bibitem{Temam1} R. Temam,  Remarks on a free boundary value problem arising in plasma physics,{\em  Commun. Partial
Differ. Equ;} 2(6), 563–585, (1977).

\bibitem{weiss1}G.S. Weiss, Partial regularity for weak solutions of an elliptic free boundary problem, {\em Commun. Partial Differ. Equ,}n 23(3–4), (1998), 439–455.
\bibitem{weiss2}G.S. Weiss, Partial regularity for a minimum problem with free boundary, {\em J. Geom. Anal,} 9(2), (1999), 317–326.

\bibitem{yang}Y. Yang, K. Perera,  Existence and nondegeneracy of ground states in critical
free boundary problems, {\em Nonlinear Anal,} 180, (2019), 75–93.
\end{thebibliography}
\end{document}